\documentclass[11pt,amssymb]{amsart}
\usepackage{amssymb}
\usepackage[mathscr]{eucal}
\usepackage[all,cmtip]{xy}

\newcommand{\Z}{{\mathbb Z}}

\newcommand{\hhe}{H^{\mathrm{\acute{e}t}}}

\newcommand{\R}{{\mathbb R}}

\newcommand{\Ga}{\mathrm{Gal}}
\newcommand{\uC}{\underline{C}}
\newtheorem{thm}{Theorem}[section]

\newcommand{\Hom}{\mathrm{Hom}}
\newcommand{\Ext}{\mathrm{Ext}}

\newcommand{\Spec}{\mathrm{Spec}}
\newcommand{\he}{H_{\mathrm{\acute{e}t}}}

\newcommand{\uCv}{\underline{C}^{(v)}}
.240pk scaled 1200 .240pk
\begin{document}

\title[Residue maps]{On residue maps for affine curves}

\author[I.~Rapinchuk]{Igor A. Rapinchuk}

\begin{abstract}
We establish several compatibility results between residue maps in \'etale and Galois cohomology that arise naturally in the analysis of smooth affine algebraic curves having good reduction over discretely valued fields. These results are needed, and in fact have already been used, for the study of finiteness properties of the unramified cohomology of function fields of affine curves over number fields.
\end{abstract}

\address{Department of Mathematics, Michigan State University, East Lansing, MI
48824, USA}

\email{rapinchu@msu.edu}

\maketitle

\section{Introduction}

The purpose of this note to establish several compatibilities between certain residue maps in \'etale and Galois cohomology that arise naturally in the analysis of affine curves with good reduction. Our main result (Theorem \ref{T-Residue}) plays a key role in \cite{CRR4}
%has been used in \cite{CRR4} as a key ingredient
in the proof of a finiteness statement for the unramified cohomology of function fields of affine curves over number fields, with applications to finiteness properties of the genus of spinor groups of quadratic forms, as well as some other groups, over such fields.
We formulate our result below in a rather general setting to make it suitable for other applications.

We begin with a description of our set-up. Suppose $k$ is a field equipped with a discrete valuation $v$, with valuation ring $\mathcal{O}_v$ and residue field $k^{(v)}$. Fix an integer $n \geq 2$ that is invertible in $\mathcal{O}_v$. Let $C$ be a smooth geometrically integral affine curve over $k$, and assume that $C$ has smooth reduction at $v$, i.e. there exists a smooth scheme $\mathcal{C}$ over $\mathcal{O}_v$ with generic fiber $C$ such that the special fiber $\uC^{(v)}$ is a smooth geometrically integral (affine) curve over $k^{(v)}.$
Then the localization sequence in \'etale cohomology for the pair $(\mathcal{C}, \uC^{(v)})$, together with absolute purity, gives a map
$$
\rho_v^{\ell} \colon H^{\ell}(C, \mu_n^{\otimes b}) \to H^{\ell-1}(\uC^{(v)}, \mu_n^{\otimes(b-1)})
$$
for all $\ell \geq 1$ and all $b \in \Z$ (see \S\ref{S-PfPartI} for the details of the construction and the end of this section for an explanation of our notations). Let us fix separable closures $\overline{k}$ and $\overline{k^{(v)}}$ of $k$ and $k^{(v)}$, respectively. Since the curves $C$ and $\uC^{(v)}$ are affine, we have $$H^q(C \otimes_k \overline{k}, \mu_n^{\otimes b}) = 0 \ \ \  \text{and} \ \ \  H^q(\uC^{(v)} \otimes_{k^{(v)}} \overline{k^{(v)}} , \mu_n^{\otimes b}) = 0$$ for $q \geq 2$ (see, e.g., \cite[Lemma 65.3]{EC-Stack}).  Consequently, the Hochschild-Serre spectral sequences
\begin{equation}\label{E-SS1}
E_2^{p,q} = H^p(k , H^q(C \otimes_k \overline{k} , \mu_n^{\otimes b})) \, \Rightarrow H^{p+q}(C , \mu_n^{\otimes b})
\end{equation}
and
\begin{equation}\label{E-SS2}
E_2^{p,q} = H^p(k^{(v)} , H^q(\uC^{(v)} \otimes_{k^{(v)}} \overline{k^{(v)}} , \mu_n^{\otimes b})) \, \Rightarrow H^{p+q}(\uC^{(v)} , \mu_n^{\otimes b})
\end{equation}
yield maps
$$
H^{\ell}(C , \mu_n^{\otimes b}) \stackrel{\omega^{b ,\ell}_k}{\longrightarrow} H^{\ell-1}(k , H^1(C \otimes_k \overline{k} , \mu_n^{\otimes b})) \ \ \ \text{and} \ \ \ H^{\ell}(\uC^{(v)} , \mu_n^{\otimes b}) \stackrel{\omega^{b ,\ell}_{k^{(v)}}}{\longrightarrow} H^{\ell-1}(k^{(v)} , H^1(\uC^{(v)} \otimes_{k^{(v)}} \overline{k^{(v)}} , \mu_n^{\otimes b}))
$$
for every $\ell \geq 1.$

Next, let us now suppose that there is an isomorphism
\begin{equation}\label{E-LocConstant}
\varphi \colon H^1(\uC^{(v)} \otimes_{k^{(v)}} \overline{k^{(v)}}, \mu_n^{\otimes b}) \stackrel{\sim}{\longrightarrow} H^1(C \otimes_k \overline{k}, \mu_n^{\otimes b})
\end{equation}
which is compatible with the actions of $G_{k^{(v)}} = \Ga(\overline{k^{(v)}}/ k^{(v)})$ and $G_k =\Ga(\overline{k}/k)$ in the following sense. Let  $I_v \subset D_v \subset G_k$ be the decomposition and inertia groups at $v$, respectively, and consider the natural isomorphism $D_v / I_v \simeq G_{k^{(v)}}$. Then we assume that $I_v$ acts
trivially on $H^1(C \otimes_k \overline{k}, \mu_n^{\otimes b})$ and $\varphi$ is a morphism of $G_{k^{(v)}}$-modules.
%i.e. such that the inertia subgroup of the decomposition group $D_v \subset G_k$ at $v$ acts trivially on $H^1(C \otimes_k \overline{k}, \mu_n^{\otimes b})$.
In particular, this means that
$H^1(C \otimes_k \overline{k}, \mu_n^{\otimes b})$ is unramified at $v$, and consequently, there exists a residue map
$$
H^{\ell}(k, H^1 (C \otimes_k \overline{k} , \mu_n^{\otimes b})) \to H^{\ell -1} (k^{(v)}, H^1 (\uC^{(v)} \otimes_{k^{(v)}} \overline{k^{(v)}}, \mu_n^{\otimes b})(-1))
$$
(see \cite[Part 1, Ch. II, \S7]{GMS}). On the other hand, in view of the well-known identification
\begin{equation}\label{E-FundGroup}
H^1 (\uC^{(v)} \otimes_{k^{(v)}} \overline{k^{(v)}}, \mu_n^{\otimes b}) \simeq \Hom (\pi_1(\uC^{(v)} \otimes_{k^{(v)}} \overline{k^{(v)}}, \bar{\eta}), \mu_n^{\otimes b}), \ \ \ \text{where} \ \bar{\eta} = {\rm Spec}~\overline{k^{(v)}(\uC^{(v)})},
\end{equation}
%(see, e.g., \cite[I.5.4]{Milne}),
we have an isomorphism $H^1 (\uC^{(v)} \otimes_{k^{(v)}} \overline{k^{(v)}}, \mu_n^{\otimes b})(-1) \simeq H^1 (\uC^{(v)} \otimes_{k^{(v)}} \overline{k^{(v)}}, \mu_n^{\otimes (b-1)})$. Thus, the above residue map can be rewritten as
\begin{equation}\label{E-ResidueGalois}
H^{\ell}(k, H^1 (C \otimes_k \overline{k} , \mu_n^{\otimes b})) \stackrel{r_v^{\ell}}{\longrightarrow} H^{\ell -1} (k^{(v)}, H^1 (\uC^{(v)} \otimes_{k^{(v)}} \overline{k^{(v)}}, \mu_n^{\otimes (b-1)})).
\end{equation}
Finally, recall that the Bloch-Ogus spectral sequence \cite{Bl-Og} %(see also \cite{CT-H-K})
yields surjective maps
$$
H^{\ell}(C, \mu_n^{\otimes b}) \stackrel{\varepsilon_C}{\longrightarrow} H^{\ell}_{ur} (k(C), \mu_n^{\otimes b}) \ \ \ \text{and} \ \ \ H^{\ell} (\uC^{(v)}, \mu_n^{\otimes b}) \stackrel{\varepsilon_{\uC^{(v)}}}{\longrightarrow} H^{\ell}_{ur} (k^{(v)}(\uC^{(v)}), \mu_n^{\otimes b}),
$$
that are induced by passage to the generic point and
where the targets are the unramified cohomology groups of the function fields with respect to the geometric places (corresponding to the closed points of the curves). Having introduced the required notations, we can now state our main result.

\begin{thm}\label{T-Residue}
Keep the notations introduced above.

\vskip1mm

\noindent {\rm (a)} \parbox[t]{17cm}{If isomorphism (\ref{E-LocConstant}) holds, then for every $\ell \geq 2$ and all integers $b$, there is a commutative diagram
\begin{equation}\tag{I}
\xymatrix{
H^{\ell}(C, \mu_n^{\otimes b}) \ar[rr]^{\omega_k^{b, \ell}} \ar[d]_{\rho_v^{\ell}} &  & H^{\ell-1}(k, H^1(C \otimes_k \overline{k} , \mu_n^{\otimes b})) \ar[d]^{\partial_v^{\ell-1}} \\ H^{\ell-1}(\uC^{(v)}, \mu_n^{\otimes(b-1)}) \ar[rr]^{\omega_{k^{(v)}}^{b-1, \ell-1}} & & H^{\ell-2}(k^{(v)}, H^1 (\uC^{(v)} \otimes_{k^{(v)}} \overline{k^{(v)}} , \mu_n^{\otimes (b-1)})),}
\end{equation}}
where $\partial_v^{\ell-1}$ coincides up to sign with $r_v^{\ell-1}.$

\vskip3mm

\noindent {\rm (b)} \parbox[t]{17cm}{For every $\ell \geq 1$ and all integers $b$, there is a commutative diagram

\begin{equation}\tag{II}
\xymatrix{
H^{\ell}(C, \mu_n^{\otimes b}) \ar[r]^{\varepsilon_C} \ar[d]_{\rho_v^{\ell}} & H^{\ell}_{ur} (k(C), \mu_n^{\otimes b}) \ar[d]^{\delta_v^{\ell}} \\ H^{\ell-1} (\uC^{(v)}, \mu_n^{\otimes(b-1)}) \ar[r]^{\varepsilon_{\uC^{(v)}}} & H^{\ell-1}_{ur} (k^{(v)}(\uC^{(v)}), \mu_n^{\otimes(b-1)}),}
\end{equation}
where $\delta_v^{\ell}$ coincides up to sign with the natural map induced by the residue map in Galois cohomology.}

\end{thm}

\vskip3mm

\noindent {\bf Remark 1.2.} We note that the isomorphism (\ref{E-LocConstant}) holds in the following situation, which was needed for the considerations in \cite{CRR4}. Let $p$ be any prime $\neq \mathrm{char}\: k^{(v)}$. Given a smooth geometrically integral affine curve $C$ over $k$, denote by $\hat{C}$ denote the smooth geometrically integral complete curve over $k$ that contains $C$ as an open subset. Let us assume that there exist models $\mathcal{C} \subset \hat{\mathcal{C}}$ of these curves over $\mathcal{O}_v$ such that the associated reductions $\uC^{(v)} \subset \underline{\hat{C}}^{(v)}$ are smooth, geometrically integral, and satisfy
$$
\vert \hat{C}(\bar{k}) \setminus C(\bar{k}) \vert = \vert \underline{\hat{C}}^{(v)}(\overline{k^{(v)}}) \setminus \uC^{(v)}(\overline{k^{(v)}}) \vert.
$$
Then the results of \cite[Ch. XIII, \S2]{SGA1} imply that the specialization map defines an isomorphism of the maximal pro-$p$ quotients of the fundamental groups
$$
\pi_1(C \otimes_k \bar{k})^{(p)} \longrightarrow \pi_1(\uC^{(v)} \otimes_{k^{(v)}} \overline{k^{(v)}})^{(p)}
$$
(with a compatible choice of base points), which yields the isomorphism (\ref{E-LocConstant}) whenever $n$ is prime to the residue characteristic $\mathrm{char}\: k^{(v)}.$
%Though not needed in the current setting, we also remark that if the model $\mathcal{C}$ over $\Spec~\mathcal{O}_v$ is smooth and {\it proper}, then the required isomorphism is an immediate consequence of the smooth and proper base change theorem (see, e.g., \cite[Ch. VI, Corollary 4.2]{Milne}).

\vskip5mm

The basic strategy for the proof of part (a) is to interpret the Hochschild-Serre spectral sequence as a Leray spectral sequence and then exploit certain formal properties of the localization sequence and the Grothendieck spectral sequence. The argument for part (b) is based on analyzing the relationship between the Kato complexes for $C$ and $\uC^{(v)}$.

The paper is organized as follows. In \S\ref{S-Derived}, we recall several facts about derived categories and spectral sequences that are needed for our arguments. Parts (a) and (b) of Theorem \ref{T-Residue} are then proved in \S\ref{S-PfPartI} and \S\ref{S-PfPartII}, respectively.

\vskip5mm

\noindent {\bf Notations and conventions.} Let $X$ be a scheme. For any positive integer
%For any field $F$ and any positive integer $n$ invertible in $F$, we denote by $\mu_n$ the group of $n$th roots of unity in a fixed separable closure $\bar{F}$ of $F$.
%Next, let $X$ be an arbitrary scheme. Then for any positive integer
$n$ invertible on $X$, we denote by $\mu_n = \mu_{n, X}$ the \'etale sheaf of $n$th roots of unity on $X$.
We follow the usual notations for the Tate twists of $\mu_n$. Namely, for $i \geq 0$, we set
$
\Z / n \Z (i) = \mu_n^{\otimes i}
$
(where $\mu_n^{\otimes i}$ is the sheaf associated to the $i$-fold tensor product of $\mu_n$), with the convention that
$
\mu_n^{\otimes 0} = \Z/n \Z.
$
If $i < 0$, we let
$$
\Z / n \Z (i) = \Hom (\mu_n^{\otimes (- i)}, \Z / n \Z).
$$
All cohomology groups considered in the paper will be \'etale cohomology. 
In the case that $X = \mathrm{Spec}~F$ for a field $F$, we identify $\mu_n$ with the group of $n$th roots of unity in a fixed separable closure $\overline{F}$ of $F$. We will also tacitly identify the \'etale cohomology of ${\rm Spec}~F$ with the Galois cohomology of $F$.

\section{Several recollections on derived categories}\label{S-Derived}

In this section, we briefly summarize several points concerning derived categories and spectral sequences that will be needed for the proof of Theorem \ref{T-Residue}(a).

First, recall that if $\mathcal{A}$ is an abelian category and $X = (X^n, d^n)_{n \in \Z}$ is a cochain complex in $\mathcal{A}$, then the truncations $\tau_{\leq n} X$ and $\tau_{\geq n} X$ are defined by the diagrams
$$
\xymatrix{\tau_{\leq n}X \ar[d]  &
\ldots \ar[r] &
X^{n - 1} \ar[r] \ar[d] &
\ker(d^n) \ar[r] \ar[d] &
0 \ar[r] \ar[d] &
\ldots \\
X  &
\ldots \ar[r] &
X^{n - 1} \ar[r] &
X^n \ar[r] &
X^{n + 1} \ar[r] &
\ldots}
$$
and
$$
\xymatrix{
X \ar[d] &
\ldots \ar[r] &
X^{n - 1} \ar[r] \ar[d] &
X^n \ar[r] \ar[d] &
X^{n + 1} \ar[r] \ar[d] &
\ldots \\
\tau_{\geq n}X &
\ldots \ar[r] &
0 \ar[r] &
{\rm coker}(d^{n - 1}) \ar[r] &
X^{n + 1} \ar[r] &
\ldots
}
$$
Furthermore, if $A$ is an object of $\mathcal{A}$, we will denote by $A^{\bullet}$ the complex that has $A$ in degree 0, and 0 everywhere else. It then follows from the definitions that for a cochain complex $X$ as above, we have $\tau_{[n]}(X) := \tau_{\geq n}(\tau_{\leq n} (X)) = H^n(X)^{\bullet}[-n].$ %We will also denote by $D^+(\mathcal{A})$ the derived category of bounded below complexes in $\mathcal{A}$.

Next, suppose
$$
\mathcal{A} \stackrel{W}{\longrightarrow} \mathcal{B} \stackrel{V}{\longrightarrow} \mathcal{C}
$$
are additive left exact functors between abelian categories with enough injectives. Recall that if $W$ maps injective objects to $V$-acyclic ones, then there is a Grothendieck spectral sequence
$$
E_2^{p,q} = R^pV(R^qW(A)) \Rightarrow R^{p+q}(VW)(A)
$$
for every bounded below complex $A$ in $\mathcal{A}.$ Furthermore, we have exact triangles
\begin{equation}\label{E-Triangle}
\tau_{\leq n} RW(A) \to RW(A) \to \tau_{\geq (n+1)} RW(A) \to \tau_{\leq n} RW(A) [1]
\end{equation}
and
\begin{equation}\label{E-Triangle1}
\tau_{\leq (n-1)} RW(A) \to \tau_{\leq n} RW(A) \to R^nW(A)[-n] \to \tau_{\leq (n-1)}RW(A)[1]
\end{equation}
in the derived category $D^+(\mathcal{B})$ of bounded below complexes in $\mathcal{B}$, and the filtration on $R^n(VW)(A)$ is given by
$$
F^j R^n(VW)(A) = {\rm Im}~\left(R^nV (\tau_{\leq (n-j)} RW(A)) \to R^nV(RW(A)) = R^n(VW)(A) \right)
$$
%where $\tau_{\leq i}$ is the trunction functor and the map comes from the natural map $\tau_{\leq (n-j)} RW(A) \to RW(A)$ in the derived category
(see, e.g., \cite[Appendix B]{Blin-Mer}).

In our situation, we consider a first quadrant spectral sequence
$$
E_2^{p,q} = R^pV(R^qW(A)) \Rightarrow H^{p+q} = R^{p+q}(VW)(A)
$$
that satisfies $R^qW(A) = 0$ for $q \geq 2$, so that, in particular, $E_2^{p,q} = 0$ for $q \geq 2$ and all $p$. For such a spectral sequence, it is clear that we have $F^{n-1}H^n = H^n$ and $E_{\infty}^{n-1,1} = E_3^{n-1,1} = \ker (d_2^{n-1,1})$ for all $n$, which yields an edge map $e \colon H^n \to E_2^{n-1,1}$ defined as the composition of the maps
$$
H^n = F^{n-1} H^n \twoheadrightarrow F^{n-1} H^n/F^n H^n = E_{\infty}^{n-1,1} \hookrightarrow E_2^{n-1,1}.
$$
On the other hand, using (\ref{E-Triangle}), it is easy to see that $R^nV(\tau_{\leq 1} RW(A)) = R^n(VW)(A)$ for all $n.$ So, since $\tau_{[1]}(A) := \tau_{\geq 1}(\tau_{\leq 1}(A)) = H^1(A)^{\bullet}[-1]$, the triangle (\ref{E-Triangle1}) allows us to conclude that $e$ is obtained by simply applying $R^nV$ to the canonical map
$$
\tau_{\leq 1} RW(A) \to \tau_{\geq 1}(\tau_{\leq 1} RW(A)).
$$

\section{Proof of Theorem \ref{T-Residue}(a)}\label{S-PfPartI}

In this section, our goal will be to establish part (a) of Theorem \ref{T-Residue} dealing with the commutativity of diagram (I).

%we turn to establishing the commutativity of diagram (I) in Theorem \ref{T-Residue}.

Let $\pi^C \colon C \to {\rm Spec}~k$ and $\pi^{\uC^{(v)}} \colon \uC^{(v)} \to {\rm Spec}~k^{(v)}$ be the structure morphisms of the curves $C$ and $\uC^{(v)}$, respectively. We first note that, as shown in \cite[Example 12.8]{Milne-LEC}, the spectral sequences (\ref{E-SS1}) and (\ref{E-SS2}) are simply the Leray spectral sequences
$$
E_2^{p,q} = H^p ({\rm Spec}~k, R^q\pi^C_*\mu_n^{\otimes d}) \Rightarrow H^{p+q}(C, \mu_n^{\otimes d})
$$
and
$$
E_2^{p,q} = H^p ({\rm Spec}~k^{(v)}, R^q \pi^{\uC^{(v)}}_* \mu_n^{\otimes d}) \Rightarrow H^{p+q}(\uC^{(v)}, \mu_n^{\otimes d}),
$$
respectively.

Next, recall that in the localization sequence for the pair $(\mathcal{C}, \uC^{(v)})$, we have the boundary map
$$
\delta_1 \colon H^{\ell}(C, \mu_n^{\otimes b}) \to H_{\uC^{(v)}}^{\ell+1} (\mathcal{C}, \mu_n^{\otimes b}).
$$
Furthermore, by absolute purity (see, e.g. \cite[Ch. VI, \S\S 5 and 6]{Milne}), there is an isomorphism
\begin{equation}\label{E-Purity1}
H_{\uC^{(v)}}^{\ell+1} (\mathcal{C}, \mu_n^{\otimes b}) \simeq H^{\ell-1}(\uC^{(v)}, \mu_n^{\otimes (b-1)}),
\end{equation}
and the map $\rho_v^{\ell} \colon H^{\ell}(C, \mu_n^{\otimes b}) \to H^{\ell-1}(\uC^{(v)}, \mu_n^{\otimes(b-1)})$ appearing on the left-hand side of diagram (I) is then defined as the composition of $\delta_1$ followed by (\ref{E-Purity1}).

To obtain a similar description for the map $\partial_v^{\ell-1}$ in (I), we first recall the following characterization of locally constant constructible (lcc) sheaves over discrete valuation rings. With $k$ and $\mathcal{O}_v$ as above, let $j \colon \Spec~k \to \Spec~\mathcal{O}_v$ denote the inclusion of the generic point. Then the functor $\mathcal{F} \mapsto j^* \mathcal{F}$ is an equivalence between the category of lcc sheaves on $\Spec~\mathcal{O}_v$ and the category of those lcc sheaves on $\Spec~k$ that correspond to finite discrete $G_k$-modules that are unramified at the closed point of $\Spec~\mathcal{O}_v$ (see, e.g., \cite[Corollary 1.1.7.3]{Conrad}). Note that the $G_k$-module corresponding to $j^*\mathcal{F}$ is simply the stalk $\mathcal{F}_{\overline{\eta}}$, where $\overline{\eta} \colon \Spec~{\overline{k}} \to \mathcal{O}_v$ is the geometric point above the generic point. Furthermore, we recall that a constructible sheaf $\mathcal{F}$ over $\Spec~O_v$ is locally constant if and only if there is an isomorphism $\mathcal{F}_{\overline{s}} \simeq \mathcal{F}_{\overline{\eta}}$ that is compatible with the respective Galois actions, where $\overline{s} \colon \Spec~\overline{k^{(v)}} \to \Spec~\mathcal{O}_v$ is the geometric point over the closed point (see \cite[V.1.10]{Milne}).

Returning to the argument, let $\pi \colon \mathcal{C} \to \Spec~\mathcal{O}_v$ be the structure morphism, and denote by $i \colon {\rm Spec}~k^{(v)} \to {\rm Spec}~\mathcal{O}_v$ and $j \colon {\rm Spec}~k \to {\rm Spec}~\mathcal{O}_v$ and by $i' \colon \uC^{(v)} \to \mathcal{C}$ and $j' \colon C \to \mathcal{C}$ the closed and open immersions for ${\rm Spec}~\mathcal{O}_v$ and $\mathcal{C}$, respectively. Consider the constructible sheaf $R^q \pi_* \mu_n^{\otimes b} \in {\rm Sh}(({\rm Spec}~\mathcal{O}_v)_{\acute{e}t})$ for $q \geq 1.$ It is easy to see (e.g. using the smooth base change theorem) that the stalk $(R^q \pi_* \mu_n^{\otimes b})_{\overline{\eta}}$, where $\overline{\eta} \colon \Spec~{\overline{k}} \to \mathcal{O}_v$ is as before, is simply the $G_k$-module $H^q(C \otimes_k \overline{k}, \mu_n^{\otimes b})$. For $q = 1$, our assumption (\ref{E-LocConstant}) guarantees that this module is unramified at the closed point; this condition holds automatically for $q \geq 2$ since in that case $H^q(C \otimes_k \overline{k}, \mu_n^{\otimes b}) = 0.$ Consequently, the preceding discussion implies that $R^q \pi_* \mu_n^{\otimes b}$ is lcc (in fact 0 for $q \geq 2$) and $(R^q \pi_* \mu_n^{\otimes b})_{\overline{s}} \simeq H^q(\uC^{(v)} \otimes_{k^{(v)}} \overline{k^{(v)}}, \mu_n^{\otimes b}).$ Therefore, we can apply absolute purity to obtain an isomorphism
\begin{equation}\label{E-Purity2}
H_{{\rm Spec}~k^{(v)}}^{\ell}({\rm Spec}~\mathcal{O}_v, R^q \pi_* \mu_n^{\otimes b}) \simeq H^{\ell-2} ({\rm Spec}~k^{(v)}, i^*(R^q \pi_* \mu_n^{\otimes b})(-1)).
\end{equation}
Note that using the usual identification between \'etale and Galois cohomology, we can rewrite the group on the right as $H^{\ell-2}(k^{(v)}, H^q (\uC^{(v)} \otimes_{k^{(v)}} \overline{k^{(v)}} , \mu_n^{\otimes b})(-1))$; for $q = 1$, this last group is simply
%and in view of (\ref{E-FundGroup}), the group on the right can be identified with
$H^{\ell -2}(k^{(v)}, H^1 (\uC^{(v)} \otimes_{k^{(v)}} \overline{k^{(v)}} , \mu_n^{\otimes (b-1)}))$ in view of (\ref{E-FundGroup}). We then define $\partial_v^{\ell-1}$ to be the composition of the boundary map
$$
\delta_2 \colon H^{\ell-1}(k, R^1\pi_* \mu_n^{\otimes b}) \to H_{{\rm Spec}~k^{(v)}}^{\ell}({\rm Spec}~\mathcal{O}_v, R^1 \pi_* \mu_n^{\otimes b})
$$
in the localization sequence of the pair $({\rm Spec}~\mathcal{O}_v, {\rm Spec}~k^{(v)})$ followed by (\ref{E-Purity2}). The results of \cite{JSS} imply that $\partial_v^{\ell-1}$ constructed in this way coincides up to sign with the residue map $r_v^{\ell-1}$ in (\ref{E-ResidueGalois}).

To continue our analysis of diagram (I), let us now briefly recall the construction of the boundary maps $\delta_1$ and $\delta_2$ (see, e.g., \cite[Theorem 9.4]{Milne-LEC} for the details). Let $\mathcal{S} = \Z_{{\rm Spec}~\mathcal{O}_v}$ be the constant \'etale sheaf on ${\rm Spec}~\mathcal{O}_{v}$ defined by $\Z.$ We have the following short exact sequence of \'etale sheaves on ${\rm Spec}~\mathcal{O}_v$
\begin{equation}\label{E-SES}
0 \to j_! j^* {\mathcal{S}} \to \mathcal{S} \to i_* i^* \mathcal{S} \to 0,
\end{equation}
where $j_!$ is the functor of extension by zero. For ease of notation, we will write $\mathcal{S}_{{\rm Spec}~k} = j_! j^* \mathcal{S}$ and $\mathcal{S}_{{\rm Spec}~k^{(v)}} = i_* i^* \mathcal{S}.$ Then for any \'etale sheaf $\mathcal{G} \in {\rm Sh}(({\rm Spec}~\mathcal{O}_v)_{\acute{e}t})$, one has canonical identifications
%$$
%H^{\ell} (C, j_C^*\mathcal{F}) \simeq {\rm Ext}^{\ell}(\Z_{{\rm Spec}~k}, \pi_* \mathcal{F}), \ \ \ \ \  H^{\ell+1}_{\uC^{(v)}} (\mathcal{C}, \mathcal{F}) \simeq {\rm Ext}^{\ell + 1}(\Z_{{\rm Spec}~k^{(v)}}, \pi_* \mathcal{F}),
%$$
%where $j_C \colon C \to \mathcal{C}$ is the open immersion, and
$$
H^{\ell} ({\rm Spec}~k, j^*\mathcal{G}) \simeq {\rm Ext}^{\ell}_{{\rm Spec}~\mathcal{O}_v}(\mathcal{S}_{{\rm Spec}~k}, \mathcal{G}) \ \ \ \text{and} \ \ \
H^{\ell}_{{\rm Spec}~k^{(v)}} ({\rm Spec}~\mathcal{O}_v, \mathcal{G}) \simeq {\rm Ext}^{\ell}_{{\rm Spec}~\mathcal{O}_v}(\mathcal{S}_{{\rm Spec}~k^{(v)}}, \mathcal{G})
$$
and $\delta_2$ is simply the boundary map in the long exact sequence induced by (\ref{E-SES}). Similarly, one has a short exact sequence of \'etale sheaves on $\mathcal{C}$
\begin{equation}\label{E-SES1}
0 \to \mathcal{T}_C \to \mathcal{T} \to \mathcal{T}_{\uC^{(v)}} \to 0,
\end{equation}
where $\mathcal{T} = \Z_{\mathcal{C}}$ and, in analogy with the preceding discussion, $\mathcal{T}_C = j'_! {j'}^* \mathcal{T}$ and $\mathcal{T}_{\uC^{(v)}} = i'_* {i'}^* \mathcal{T}.$ As above, for any \'etale sheaf $\mathcal{F} \in {\rm Sh}(\mathcal{C}_{\acute{e}t})$, we have identifications
$$
H^{\ell}(C, {j'}^*\mathcal{F}) \simeq {\rm Ext}^{\ell}_{\mathcal{C}} (\mathcal{T}_{C}, \mathcal{F}) \ \ \ \text{and} \ \ \ H^{\ell}_{\uC^{(v)}} (\mathcal{C}, \mathcal{F}) \simeq \Ext^{\ell}_{\mathcal{C}} (\mathcal{T}_{\uC^{(v)}}, \mathcal{F}),
$$
with $\delta_1$ arising as the boundary map induced by (\ref{E-SES1}). We also note the isomorphisms
$$
\Ext_{\mathcal{C}}^{\ell}(\mathcal{T}_C, \mathcal{F}) \simeq \Ext_{\mathcal{C}}^{\ell}(\pi^*\mathcal{S}_{{\rm Spec}~k}, \mathcal{F}) \ \ \ \text{and} \ \ \ \Ext^{\ell}_{\mathcal{C}} (\mathcal{T}_{\uC^{(v)}}, \mathcal{F}) \simeq \Ext^{\ell}_{\mathcal{C}} (\pi^*\mathcal{S}_{{\rm Spec}~k^{(v)}} , \mathcal{F}).
$$
Thus, to complete the proof of Theorem \ref{T-Residue}(a), it suffices to show that we have the following commutative diagram
\begin{equation}\label{E-Diagram}
\xymatrix{\Ext_{\mathcal{C}}^{\ell}(\pi^* \mathcal{S}_{{\rm Spec}~k}, \mu_n^{\otimes b}) \ar[d] \ar[r] & \Ext_{{\rm Spec}~\mathcal{O}_v}^{\ell-1}(\mathcal{S}_{{\rm Spec}~k}, R^1 \pi_* \mu_n^{\otimes b}) \ar[d] \\ \Ext^{\ell +1}_{\mathcal{C}} (\pi^* \mathcal{S}_{{\rm Spec}~k^{(v)}}, \mu_n^{\otimes b}) \ar[r] & \Ext_{{\rm Spec}~\mathcal{O}_v}^{\ell}(\mathcal{S}_{{\rm Spec}~k^{(v)}}, R^1 \pi_* \mu_n^{\otimes b})}
\end{equation}
where the vertical maps are the boundary maps discussed above, and the horizontal maps are edge maps in the appropriate Leray spectral sequences.

For this, it will be convenient to pass to the framework of derived categories. Denote by $D^+(\mathcal{C})$, $D^+({\rm Spec}~\mathcal{O}_v)$, and $D^+(Ab)$ the derived categories of bounded below complexes of \'etale sheaves on $\mathcal{C}$, \'etale sheaves on ${\rm Spec}~\mathcal{O}_v$, and abelian groups, and consider the derived functors
$$
{\bf R}\pi_* \colon D^+(\mathcal{C}) \to D^+({\rm Spec}~\mathcal{O}_v) \ \ \ \text{and} \ \ \ {\bf R}\Hom_{{\rm Spec}~\mathcal{O}_v} (\mathcal{F}, \cdot) \colon D^+({\rm Spec}~\mathcal{O}_v) \to D^+ (Ab),
$$
where $\mathcal{F}$ is an \'etale sheaf on ${\rm Spec}~\mathcal{O}_v$.
%, and use the fact that  %${\bf R} \pi_*$
%$({\bf R}\Hom_{{\rm Spec}~\mathcal{O}_v} (\mathcal{F}, \cdot)) \circ ({\bf R} \pi_*)$ coincides with the derived functor ${\bf R}\Hom_{\mathcal{C}}(\pi^* \mathcal{F}, \cdot) \colon D^+(\mathcal{C}) \to D^+(Ab).$
Now, starting with the
commutative diagram
\begin{equation}\label{E-Diagram1}
\xymatrix{{\bf R}\Hom_{\Spec~\mathcal{O}_v}(\mathcal{S}_{\Spec~k}, \tau_{\leq 1}({\bf R}\pi_* \mu_n^{\otimes b})) \ar[r] \ar[d] & {\bf R}\Hom_{\Spec~\mathcal{O}_v}(\mathcal{S}_{\Spec~k}, \tau_{[1]}({\bf R}\pi_* \mu_n^{\otimes b})) \ar[d] \\ {\bf R}\Hom_{\Spec~\mathcal{O}_v}(\mathcal{S}_{\Spec~k^{(v)}}, \tau_{\leq 1}({\bf R}\pi_* \mu_n^{\otimes b}))[1] \ar[r] & {\bf R}\Hom_{\Spec~\mathcal{O}_v}(\mathcal{S}_{\Spec~k^{(v)}}, \tau_{[1]}({\bf R}\pi_* \mu_n^{\otimes b}))[1]}
\end{equation}
in $D^+(Ab)$, we take cohomology in degree $\ell$ to obtain the following commutative diagram of abelian groups
\begin{equation}\label{E-Diagram2}
\xymatrix{\R^{\ell}\Hom_{\Spec~\mathcal{O}_v}(\mathcal{S}_{\Spec~k}, \tau_{\leq 1}({\bf R}\pi_* \mu_n^{\otimes b})) \ar[r] \ar[d] & \R^{\ell}\Hom_{\Spec~\mathcal{O}_v}(\mathcal{S}_{\Spec~k}, \tau_{[1]}({\bf R}\pi_* \mu_n^{\otimes b})) \ar[d] \\ \R^{\ell}(\Hom_{\Spec~\mathcal{O}_v}(\mathcal{S}_{\Spec~k^{(v)}}, \tau_{\leq 1}({\bf R}\pi_* \mu_n^{\otimes b}))[1]) \ar[r] & \R^{\ell}(\Hom_{\Spec~\mathcal{O}_v}(\mathcal{S}_{\Spec~k^{(v)}}, \tau_{[1]}({\bf R}\pi_* \mu_n^{\otimes b}))[1])}
\end{equation}
We claim this is precisely the required diagram (\ref{E-Diagram}). Indeed, the above discussion shows that
$$
\R^{\ell}\Hom_{\Spec~\mathcal{O}_v}(\mathcal{S}_{\Spec~k}, \tau_{\leq 1}({\bf R}\pi_* \mu_n^{\otimes b})) = \R^{\ell}\Hom_{\Spec~\mathcal{O}_v}(\mathcal{S}_{\Spec~k}, ({\bf R}\pi_* \mu_n^{\otimes b})) = \Ext_{\mathcal{C}}^{\ell}(\pi^* \mathcal{S}_{{\rm Spec}~k}, \mu_n^{\otimes b}),
$$
where the last equality follows from the fact that $({\bf R}\Hom_{{\rm Spec}~\mathcal{O}_v} (\mathcal{S}_{\Spec~k}, \cdot)) \circ ({\bf R} \pi_*)$ coincides with ${\bf R}\Hom_{\mathcal{C}}(\pi^* \mathcal{S}_{\Spec~k}, \cdot).$ By the same argument,
$$
\R^{\ell}(\Hom_{\Spec~\mathcal{O}_v}(\mathcal{S}_{\Spec~k^{(v)}}, \tau_{\leq 1}({\bf R}\pi_* \mu_n^{\otimes b}))[1]) = \Ext^{\ell +1}_{\mathcal{C}} (\pi^* \mathcal{S}_{{\rm Spec}~k^{(v)}}, \mu_n^{\otimes b}).
$$
Furthermore, we have $\tau_{[1]}({\bf R}\pi_* \mu_n^{\otimes b}) = (R^1 \pi_* \mu_n^{\otimes b})^{\bullet}[-1]$, from which it follows that
$$
\R^{\ell}\Hom_{\Spec~\mathcal{O}_v}(\mathcal{S}_{\Spec~k}, \tau_{[1]}({\bf R}\pi_* \mu_n^{\otimes b})) = \Ext_{{\rm Spec}~\mathcal{O}_v}^{\ell-1}(\mathcal{S}_{{\rm Spec}~k}, R^1 \pi_* \mu_n^{\otimes b})
$$
and similarly
$$
\R^{\ell}(\Hom_{\Spec~\mathcal{O}_v}(\mathcal{S}_{\Spec~k^{(v)}}, \tau_{[1]}({\bf R}\pi_* \mu_n^{\otimes b}))[1]) = \Ext_{{\rm Spec}~\mathcal{O}_v}^{\ell}(\mathcal{S}_{{\rm Spec}~k^{(v)}}, R^1 \pi_* \mu_n^{\otimes b}).
$$
Thus, all of the terms in diagrams (\ref{E-Diagram}) and (\ref{E-Diagram2}) match up, and the vertical maps clearly coincide. Moreover, the description of the edge maps recalled in \S\ref{S-Derived} implies that the top and bottom horizontal maps in (\ref{E-Diagram2}) are precisely the edge maps in the corresponding Leray spectral sequences, as needed. This completes the proof.

\section{Proof of Theorem \ref{T-Residue}(b)}\label{S-PfPartII}

We now turn to the proof part (b) of Theorem \ref{T-Residue}. The main point of the argument is to consider a certain residue map between the Kato complexes of $C$ and $\uC^{(v)}.$ Such a map was defined by Kato in his work on cohomological Hasse principles (see \cite[\S 5]{Kato}), and then was treated more systematically using Bloch-Ogus spectral sequences by Jannsen and Saito \cite{JS}. We note that although all arguments can be carried out in the setting of \'etale and Galois cohomology (as in \cite{CT-H-K}), due to the homological nature of the Kato complex, it is somewhat more convenient to work with \'etale homology following the approach of Jannsen and Saito, which we briefly recall. To keep our notations consistent with those of \cite{JS}, we will write $\Z/ n \Z (b)$ for $\mu_n^{\otimes b}.$

Suppose $\mathfrak{C}$ is a category of noetherian schemes such that for any object $X$ of $\mathfrak{C}$, all closed immersions $i \colon Y \hookrightarrow X$ and all open immersions $j \colon V \hookrightarrow X$ are morphisms in $\mathfrak{C}.$ Let $\mathfrak{C}_*$ be the category having the same objects as $\mathfrak{C}$, but where the morphisms are only the proper maps in $\mathfrak{C}.$ We recall that a {\it homology theory} on $\mathfrak{C}$ (in the sense of Bloch-Ogus) is a sequence of covariant functors
$$
H_a (\cdot) \colon \mathfrak{C}_* \to \text{abelian groups} \ \ \ \ \ \ (a \in \Z)
$$
satisfying

\vskip1mm

\noindent $\bullet$ \parbox[t]{16cm}{{\it Contravariance under open immersions}: for any open immersion $j \colon V \hookrightarrow X$ in $\mathfrak{C}$, there is a functorial map $j^* \colon H_a (X) \to H_a (V)$; and}

\vskip1mm

\noindent $\bullet$ \parbox[t]{16cm}{{\it Localization}: if $i \colon Y \hookrightarrow X$ is a closed immersion with open complement $j \colon V \hookrightarrow X$, there is a long exact localization sequence
$$
\cdots \stackrel{\delta}{\longrightarrow} H_a (Y) \stackrel{i_*}{\longrightarrow} H_a (X) \stackrel{j^*}{\longrightarrow} H_a (V) \stackrel{\delta}{\longrightarrow} H_{a-1} (Y) \to \cdots
$$
that is functorial with respect to proper maps and open immersions.}

\vskip3mm

\noindent {\bf Examples 4.1.} Here are several constructions that are relevant for our discussion; the reader is referred to \cite[2.2, 2.3, and 2.5]{JS} for the details.

\vskip1mm

\noindent (a) Let $S$ be a noetherian scheme and $Sch_{sft/S}$ the category of separated schemes of finite type over $S$. Given a bounded complex $\Lambda$ of \'etale sheaves on $S$, one obtains a homology theory on $Sch_{sft/S}$ (the so-called {\it \'etale homology}) by setting
$$
\hhe_a (X/S, \Lambda) = \he^{-a} (X, Rf^! \Lambda),
$$
where $Rf^! \colon D^+ (S_{\text{\'et}}) \to D^+ (X_{\text{\'et}})$ is the exceptional inverse image functor (see \cite[Expos\'e XVIII, 3.1.4]{SGA4}).

\vskip1mm

\noindent (b) If $H$ is a homology theory, then for any integer $N$, one can define a shifted homology theory $H[N]$ by setting $H[N]_a (Z) = H_{a+N} (Z)$ and multiplying the connecting morphisms $\delta$ in the localization sequence by $(-1)^N.$

\vskip1mm

\noindent (c) For any scheme $X$ in $\mathfrak{C}$, let $\mathfrak{C}/X$ be the subcategory of schemes over $X$. If $H$ is a homology theory on $\mathfrak{C}/X$, then for any immersion $Z \hookrightarrow X$, we have a homology theory $H^{(Z)}$ on $\mathfrak{C}/X$ defined by setting
$$
H_a^{(Z)} (T) = H_a (T \times_X Z).
$$
Moreover, if $Y$ is a closed subscheme of a noetherian scheme $X$ with open complement $U = X \setminus Y$, then  there is a morphism of homology theories
\begin{equation}\label{E-HomologyTheory}
\delta \colon H^{(U)}[1] \to H^{(Y)}
\end{equation}
induced by the connecting morphism $\delta \colon H_a (T \times_X U) \to H_{a-1}(T \times_X Y)$ for $T \in \mathfrak{C}/X$. %(cf. \cite[Corollary 2.5]{JS}).

\vskip4mm

An important observation, due to Bloch and Ogus \cite{Bl-Og}, is that there is a {\it niveau spectral sequence} associated to any homology theory: more precisely, if $H$ is a homology theory on $\mathfrak{C}$, then for every object $X$ of $\mathfrak{C}$, there is a spectral sequence of homological type
\begin{equation}\label{E-NiveauSS}
E^1_{r,q} (X) = \bigoplus_{x \in X_r} H_{r+q}(x) \Rightarrow H_{r+q}(X),
\end{equation}
where $X_r = \{ x \in X \mid \dim \overline{\{ x \}} = r \}$ and
$$
H_a (x) = \lim_{\rightarrow} H_a (V),
$$
with the limit taken over all open non-empty subschemes $V \subset \overline{ \{x  \}}.$ In the case where our homology theory is \'etale homology, we will denote by
$$
E^1_{r,q} (X/S, \Z/n\Z(b)) = \bigoplus_{x \in X_r} \hhe_{r+q} (x/S, \Z/n\Z(b)) \Rightarrow \hhe_{r+q}(X/S, \Z/ n \Z(b))
$$
the associated niveau spectral sequence for a scheme $f \colon X \to S$ in $Sch_{sft/S}$ and the \'etale sheaf $\Z / n \Z(b)_S$ on $S$.
Now, using the construction of (\ref{E-NiveauSS}) via exact couples arising from appropriate localization sequences, Jannsen and Saito noted that if $S = \Spec~(A)$ for a discrete valuation $A$, and $X$ is a separated scheme of finite type over $S$ with generic fiber $X_{\eta}$ and special fiber $X_{s}$, then for any homology theory $H$ on the category $Sub(X)$ of subschemes of $X$ (considered as schemes over $X$), the morphism $\delta \colon H^{(X_{\eta})}[1] \to H^{(X_s)}$ from (\ref{E-HomologyTheory}) induces a morphism of spectral sequences
\begin{equation}\label{E-SpecSeqMap}
\Delta_X \colon E^1_{r,q}(X_{\eta})^{(-)} \to E^1_{r, q-1}(X_s),
\end{equation}
where the superscript $(-)$ means that all differentials in the original spectral sequence are multiplied by -1 (see \cite[Proposition 2.12]{JS}).

Furthermore, in the case where $X$ is a smooth geometrically integral variety of dimension $d$ over a field $F$ and $n$ is invertible in $F$, we recall that Poincar\'e duality yields canonical isomorphisms
$$
\hhe_a(x/F, \Z / n \Z (b)) \simeq H^{2r - a}(k(x), \Z/ n \Z (r-b)) \ \ \ \text{for} \ x \in X_r,
$$
where $k(x)$ is the residue field of $x$, and
$$
\hhe_a (X/F, \Z / n \Z (b) \simeq \he^{2d - a}(X, \Z / n \Z(d-b))
$$
(cf. \cite[Theorem 2.14]{JS}). Thus, from the complex $E^1_{\bullet, q}(X/F, \Z / n \Z (b))$, we deduce the complex
$$
\cdots \to \bigoplus_{x \in X_r} H^{r-q} (k(x), \Z / n \Z (r-b)) \to \bigoplus_{x \in X_{r-1}} H^{r-q-1} (k(x), \Z/ n \Z (r-b-1)) \to \cdots
$$
$$
\cdots \to \bigoplus_{x \in X_0} H^{-q} (k(x), \Z / n \Z (-b))
$$
which, as shown in \cite[Theorem 2.5.10]{JSS}, coincides up to sign with Kato's complex constructed using residue maps in Galois cohomology. In particular, it follows that
\begin{equation}\label{E-UnramBlochOgus}
E^2_{d,q}(X/F, \Z / n \Z(b)) = H^{d-q}_{ur}(F(X), \Z/n\Z(d-b)).
\end{equation}

With these preliminaries, let us now return to the proof of Theorem \ref{T-Residue}(b). Recall that we consider a smooth scheme $\mathcal{C}$ over $\Spec~\mathcal{O}_v$ with generic fiber $C$ and special fiber $\uC^{(v)}.$ First, viewing $C$ and $\uC^{(v)}$ as schemes over $\Spec~\mathcal{O}_v$ on the one hand, and as curves over $\Spec~k$ and $\Spec~k^{(v)}$, respectively, on the other, we have the following isomorphisms of niveau spectral sequences
$$
E^1_{r,q}(C/\Spec~\mathcal{O}_v, \Z/ n \Z(b)) \simeq E^1_{r,q}(C/k, \Z / n \Z(b))
$$
and
$$
E^1_{r,q}(\uC^{(v)}/\Spec~\mathcal{O}_v, \Z / n \Z (b)) \simeq E^1_{r, q+2}(\uC^{(v)}/k^{(v)}, \Z /n \Z(b+1))
$$
(see \cite[Lemma 2.17]{JS}). Combining this with (\ref{E-SpecSeqMap}), we thus obtain a morphism
\begin{equation}\label{E-SpecSeqResidue}
\Delta_{\mathcal{C}} \colon E^1_{r,q}(C/k, \Z / n \Z(b))^{(-)} \to E^1_{r,q+1}(\uC^{(v)}/k^{(v)}, \Z / n \Z(b+1))
\end{equation}
of spectral sequences. So, taking $q = 1 - \ell$ and replacing $b$ by $1-b$, the morphism (\ref{E-SpecSeqResidue}), in conjunction with the isomorphism (\ref{E-UnramBlochOgus}), yields a map
$$
\delta^{\ell}_v \colon H^{\ell}_{ur}(k(C), \Z/ n \Z (b)) \to H^{\ell-1}_{ur}(k^{(v)}(\uCv), \Z / n \Z (b-1)),
$$
which, according to \cite[Lemma 2.20]{JS}, coincides up to sign with the usual residue map in Galois cohomology. Also, since
$E^2_{r,q}(C/k, \Z / n \Z(b)) = 0$ for $r \neq 0,1,$ there is an edge map
\begin{equation}\label{E-Edge}
\varepsilon_C \colon H_t(C/k, \Z / n \Z (b)) \to E^2_{1, t-1}(C/k, \Z / n \Z (b)),
\end{equation}
for any $t$; in view of the isomorphism between \'etale homology and cohomology recalled above, for $t = 2- \ell$ we thus obtain a map
$$
\varepsilon_C \colon \he^{\ell}(C, \Z / n \Z (b)) \to H^{\ell}_{ur}(k(C), \Z/ n \Z (b)).
$$
Similarly, we have a map
$$
\varepsilon_{\uCv} \colon H^{\ell-1}(\uCv, \Z / n \Z (b-1)) \to H^{\ell-1}_{ur}(k^{(v)}(\uCv), \Z / n \Z (b-1)).
$$
We also note that the morphism of spectral sequences (\ref{E-SpecSeqResidue}) yields a map
$$
\rho_v^{\ell} \colon H^{\ell}(C, \Z / n \Z (b)) \to H^{\ell-1}(\uCv, \Z / n \Z (b-1))$$
of limiting terms that fits into the commutative diagram
$$
\xymatrix{H^{\ell}(C, \Z/ n \Z(b)) \ar[r]^{\varepsilon_C} \ar[d]_{\rho_v^{\ell}} & H^{\ell}_{ur} (k(C), \Z / n \Z(b)) \ar[d]^{\delta_v^{\ell}} \\ H^{\ell-1} (\uCv, \Z / n \Z(b-1)) \ar[r]^{\varepsilon_{\uCv}} & H^{\ell-1}_{ur} (k^{(v)}(\uCv), \Z \ n \Z(b-1))}
$$
Finally, it remains to note that the construction of the niveau spectral sequence implies that both $\varepsilon_C$ and $\varepsilon_{\uCv}$ are induced by passage to the generic point, whereas $\rho_v^{\ell}$ coincides with the map derived from the localization sequence that was discussed in \S\ref{S-PfPartI}. This concludes the proof. $\Box$

\vskip2mm

\noindent {\small  {\bf Acknowledgements.} I would like to thank A.S.~Merkurjev for very helpful discussions on derived categories. I would also like to thank V.I.~Chernousov, G.~Pappas, and C.~Weibel for useful comments and correspondence. I was partially supported by an AMS-Simons Travel Grant.}

\vskip5mm

\bibliographystyle{amsplain}

\end{document}